\documentclass[11pt]{article}
\usepackage[margin=0.6in]{geometry}
\usepackage{amsmath, amssymb,amsthm}  
\usepackage{listings}   
\usepackage{enumitem}
\usepackage{booktabs,bm}  
\usepackage{graphicx}   
\usepackage[dvipsnames]{xcolor}        
\usepackage{array}
\usepackage{comment}
\usepackage{lscape, wrapfig}

\usepackage[final]{pdfpages}
\newcommand{\sgn}{\operatorname{sgn}}

\title{Two intriguing variants of the AAA algorithm for rational approximation}
\author{Will Mitchell}
\begin{document}
\pagenumbering{gobble}
\maketitle

\abstract{We consider the problem of finding a rational function in barycentric form to approximate a given function or data set in $\mathbb{R}$ or $\mathbb{C}$. The famous AAA algorithm, introduced in 2018, constructs such a rational function: the barycentric weights are the entries of the final right singular vector of a \emph{Loewner matrix}, which has more rows than columns.  
We present two variants of the AAA algorithm, inspired by two intriguing quotations from the original paper. In the first, which we call {\tt AAAsmooth}, we take the barycentric weights to be a complex linear combination of the last \emph{two} right singular vectors, which eliminates the problem of spurious poles in real-valued problems and yields smoother convergence curves.  
In the second, {\tt AAAbudget}, we incorporate first derivative information. 
This allows us to use a smaller, square alternative to the Loewner matrix, so the SVDs are cheaper while the resulting approximant is similar to the result of standard AAA. We present numerical tests showing that while both variants behave fairly similarly to standard AAA, {\tt AAAsmooth} can give somewhat better results and {\tt AAAbudget} can be much faster. 
\section{Introduction}
\subsection{Rational functions in barycentric form} 
A \emph{barycentric formula} for a rational function $r(z)$, with nodes $z_j$, weights $w_j$, and data $f_j$, is given by  
\begin{equation}
r(z) = \left(\sum_{j=1}^N \frac{w_jf_j}{z-z_j}\right)\Big/\left(\sum_{j=1}^N \frac{w_j}{z-z_j}\right).
\label{eq:r(z)}
\end{equation}
Thus $r(z) = n(z)/d(z)$ where $n$ and $d$ are rational functions given in partial-fractions form.  
We assume $N\ge 1$ and we allow the nodes, weights and data to be arbitrary complex numbers, except that the nodes must be distinct and none of the weights may vanish. The nodes are also called \emph{support points}. Note that both $n$ and $d$ have simple poles at the support points, with residues $w_jf_j$ and $w_j$, respectively.
Typically the $f_j$ are known values of a given function that we seek to approximate, that is, $f_j = f(z_j)$ and we want $r(z)\approx f(z)$. 
We emphasize some important facts about $r(z)$: 

\begin{enumerate}
\item \emph{Constant data:} If we have $f_j\equiv F$ for all $j$, then $F$ can be factored from the numerator and we find $r(z)= F$ for all $z$. The weights do not matter in this case.  
\item \emph{Normalization of weights:} There is no change in $r(z)$ if we multiply all of the weights by a constant. Therefore, we may assume without loss of generality that the weights satisfy a normalization condition such as $\|w\|_2 = 1$ or $\|w\|_\infty = 1$. 
\item \emph{Exact interpolative property:} There is a removable discontinuity at each node. Indeed, if $z\approx z_j$ then one term dominates the sum in the numerator, and the same term without the factor of $f_j$ dominates the sum in the denominator. Thus the quotient $r(z)$ satisfies $\lim_{z\to z_j}r(z) = f_j$. In particular, if we use the barycentric formula \eqref{eq:r(z)} as the definition of $r(z)$ for $z\neq z_j$ and then define $r(z_j)=f_j$, the function $r(z)$ is analytic at each $z_j$. 
\item \emph{Numerical interpolative property:} It is safe to evaluate $r(z)$ for $z\approx z_j$ using equation \eqref{eq:r(z)}, even in the presence of rounding errors. For example, suppose that $z$ and $z_j$ are distinct floating-point numbers and $|z-z_j|$ is on the order of machine precision. In this case the computed value of $1/(z-z_j)$ might not have any correct digits, but it will be large, and -- crucially -- it appears in both the numerator $n(z)$ and the denominator $d(z)$. The magnitude of this quantity makes the other terms in $n(z)$ and $d(z)$ insignificant, and then the imprecise values cancel each other in the division $n/d$, yielding $f_j$ as desired. This is a key advantage of barycentric representations of rational functions.\footnote{On the other hand, \eqref{eq:r(z)} is vulnerable to loss of precision near the zeros of $n(z)$ and particularly near the zeros of $d(z)$ since these zeros are the result of cancellations in the sums over $j$. More formally, see \cite{higham2004numerical} for a proof of forward stability.}  
\item \emph{Generality:} Any rational function which takes finite values $f_j$ at the nodes $z_j$, and which is equal to a ratio of polynomials of degree at most $N-1$, can be expressed in the form \eqref{eq:r(z)} through an appropriate choice of weights $w_j$. This is equivalent to Theorem 2.1 in \cite{nakatsukasa2018aaa}. To appreciate the value of this statement, note that the simple partial-fractions expression $\sum_j w_j / (z - z_j)$ is incapable of representing arbitrary rational functions: for example, no parameter choice can yield $1/(z-2)^2$ exactly since the pole at $z=2$ is second-order. In contrast, we can readily find barycentric representations of $1/(z-2)^2$. One example using nodes $0$, $1$, and $3$ is:
\begin{equation}
\frac{1}{(z-2)^2} = \left(\frac{0.25\cdot 8}{z-0}+\frac{1\cdot(-3)}{z-1} + \frac{1\cdot 1}{z-3}\right)\Big/\left(\frac{8}{z-0}+\frac{-3}{z-1} + \frac{1}{z-3}\right)
\end{equation}
\item \emph{Polynomials in barycentric form:} Since polynomials are rational functions, the preceding statement implies that it is possible to choose the weights in such a way that $r(z)$ is the unique \emph{polynomial} interpolant of the data $f_j$ at locations $z_j$. In general, the weights for polynomial interpolation are given by the formula $w_j = \prod_{k\neq j} 1/(z_j-z_k)$. Closed formulas for the weights are known in some special cases, e.g. if the nodes are evenly spaced, or if the nodes form a Chebyshev grid or a Gauss-Legendre quadrature grid on $[-1,1]$; see \cite{doi:10.1137/S0036144502417715,trefethen2019approximation}. 
\item  \emph{First derivative at nodes:} At a support point $z_i$, the value of the first derivative of $r$ is given by 
\begin{equation}
r'(z_i) = \frac{-1}{w_i}\sum_{j\neq i}\frac{f_i-f_j}{z_i-z_j}w_j
\label{eq:r'(z)}
\end{equation}
This is a specialization of an arbitrary-order derivative formula due to Schneider and Warner \cite{schneider1986some}. Here we do not need higher-order derivatives, nor do we need to evaluate the first derivative away from the support points. 
\item \emph{Froissart doublets}: It can happen that both $n(z)$ and $d(z)$ vanish at some $z=Z$. In this case we say that $r$ has a Froissart doublet at $Z$. If we express $r$ as a ratio of polynomials of degree $N-1$, they will have a common root and so it should be possible to find another barycentric representation of this rational function using one fewer barycentric node. 
Numerically, it might instead occur that $n$ and $d$ have nearly-coinciding zeros, so that $r$ has a pole of small residue very close to a zero. Evaluation of \eqref{eq:r(z)} near a Froissart doublet may lead to some loss of precision in machine arithmetic because of cancellation error; an imperfect alternative is to use l'H\^opital's rule, $r(Z) = n'(Z)/d'(Z)$. 
For these reasons, Froissart doublets are usually viewed as undesirable. 
\item \emph{Zeros, poles and residues}: The zeros and poles of $r(z)$ can be determined from two generalized $N\times N$ eigenvalue problems. Once the poles are known, the residues can be determined by contour integration around small circular paths centered at the poles. For details, see \cite{nakatsukasa2018aaa}. 
\end{enumerate}

\subsection{The standard AAA algorithm}
Here we summarize the usual AAA procedure, emphasizing the features that we will later modify. Readers who are new to the AAA algorithm should also see the original presentation \cite{nakatsukasa2018aaa}. Our goal is to construct an approximant $r(z)$, in barycentric form, to a given target function $f(z)$. More precisely, we are given a tolerance $\epsilon$, typically on the order of $10^{-13}$, and a finite table of values of $f$ at some points in $\mathbb{C}$, and we want to choose the degree $N$, the nodes $z_j$, and the weights $w_j$ such that $r$ given by \eqref{eq:r(z)} satisfies $|r(z)-f(z)|<\epsilon$ for all $z$ within our data set. The AAA algorithm is an iterative procedure where, at each stage, $N$ is increased by one, a new node $z_j$ is introduced, and the weights $w_j$ are determined through numerical linear algebra; all of the weights are adjusted at each step, but the nodes never move. The algorithm starts with two nodes: $z_1$ is the location where $f$ differs the most from its mean value, and $z_2$ is the location where $f$ differs the most from $f(z_1)$.

The step where we add a new node is a matter of greed. Suppose that we have constructed a barycentric representation $r(z)$ given by \eqref{eq:r(z)}. If the condition $|r-f|<\epsilon$ is not satisfied, we find the location where $|r-f|$ is greatest and we let this become a new barycentric node. In the standard formulation this is not a continuous maximization problem, but rather a discrete one, since our information about $f$ comes in the form of a finite list of known function values.  

In the other step, we seek barycentric weights. 
Note that for any choice of nonzero weights $w_j$, the interpolative property ensures that $r(z_j)=f(z_j)$. We therefore need additional information about $f$ to determine the weights. In the original AAA algorithm this additional information consists of function values at other locations, called \emph{sample points}; these are the locations in our data set which have not yet been selected as support points. Suppose that $Z$ is a sample point, so $Z$ is distinct from all $z_j$. The equation $r(Z)=f(Z)$ is nonlinear in $w$, but we can instead seek to enforce the linear equation $f(Z)d(Z)=n(Z)$, that is,  
\begin{equation}
f(Z)\sum_{j=1}^N \frac{w_j}{Z-z_j} = \sum_{j=1}^N \frac{w_jf_j}{Z-z_j}
\end{equation}
which can be rearranged to obtain 
\begin{equation}
\sum_{j=1}^N \frac{f(Z) - f_j}{Z-z_j}w_j  = 0.
\end{equation}
Thus the vector of weights is perpendicular to the vector with entries $(f(Z)-f_j)/(Z-z_j)$. Note that these entries are first order finite-difference approximations of $f'$; for purely real problems they are the slopes of the secant lines that join $(Z,f(Z))$ and $(z_j, f_j)$ on the graph of $f$. 

Given data at many sample points, say $f(Z_i)=F_i$ for $i\le M$, we assemble the $M\times N$ \emph{Loewner matrix} $A$ whose entries are $A_{ij} = (F_i-f_j)/(Z_i-z_j)$. The vector of weights is perpendicular to each row of $A$, that is, the weights vector lies in the nullspace of $A$. Typically the sample set is larger than the support set, so $A$ has more rows than columns and there is no nullspace, at least in exact arithmetic. In practice, $A$ may be ill-conditioned and we can hope to find a weights vector $w$, without any vanishing entries, in the numerical nullspace of $A$. This guarantees that $r(z)$ will take the same values as $f$ at the sample points (or nearly so) as well as matching $f$ perfectly at the support points. Even if $A$ has full numerical rank, we can minimize errors in $n\approx fd$ by finding the (reduced) singular value decomposition and letting the weights be the entries of the right singular vector which corresponds to the smallest (last) singular value. That is, we let $A = U\Sigma V^H$ and then set the weights vector equal to the last column of $V$. We can then reassess whether $|r-f|<\epsilon$ and continue by adding another barycentric node if the tolerance has not been reached. 

The iteration terminates when $|r-f|<\epsilon$ at all locations where we have information about $f$. In this paper we do not use any cleanup or postprocessing procedure, e.g. for removal of Froissart doublets. 

The most expensive step in the AAA procedure is the singular value decomposition of the $M\times N$ Loewner matrix, which must be recomputed at every step. Note that $M$ decreases and $N$ increases as the algorithm progresses, because a sample point becomes a support point at each step. For many problems the algorithm terminates while $N$ is still fairly small (in the dozens, not the hundreds or thousands), so AAA usually converges quickly, often in less than a second. 

We note two undesirable ways in which the AAA process can terminate. First, if only a small number of samples are given or the tolerance is too tight, it is possible that the stopping criterion will not be satisfied and the Loewner matrix will get shorter and wider until $M<N$. If this occurs, we use a full (not reduced) SVD of $A$ to get a nullspace vector and we are guaranteed an approximant with $n(Z_i)=f(Z_i)d(Z_i)$ at all sample points; sadly this approximant risks overfitting the data. This can be avoided by taking a larger number of samples, loosening the tolerance, or imposing a maximum iteration ceiling. Second, it is possible for the last right singular vector of $A$ to contain a zero entry, resulting in a vanishing barycentric weight, say $w_k=0$. As a result, the corresponding terms in $n(z)$ and $d(z)$ vanish and we have an approximant of lower degree. Lamentably, we no longer have a guarantee that $|f_k-r(z_k)|<\epsilon$. Indeed, $z_k$ is not a support point (since its weight vanished) so we lack  the interpolative property $r(z_k)=f_k$, and simultaneously the Loewner matrix does not contain any row corresponding to $z_k$, so we are not enforcing the orthogonality condition $n(z_k)\approx f_kd(z_k)$ either. This danger seems to have not yet been remarked on in the literature, possibly reflecting its rarity in practice. In our tests, zero weights do sometimes appear in the last step of the AAA process, especially in cases where the convergence has stagnated and become noisy before reaching the tolerance; however, it is extremely rare to actually meet the convergence threshold at all other data points while failing at $z_k$. We have not seen any example with $|f_k-r(z_k)| > \epsilon > \max_{j\neq k} |f_j-r(z_j)|$, although the approximation of $\tan(256z)$ in Fig. \ref{fig:tan256} below comes close.  

Since 2018 the AAA algorithm has stimulated a large and diverse body of generalizations and applications. To conclude this summary of the AAA algorithm, we mention only two subsequent works: the trigonometric version for periodic problems \cite{baddoo2021aaatrig} and the continuum version \cite{driscoll2024aaa}, where the sample points are adaptively placed in between adjacent support points, allowing for automatic refinement of the grid near challenging features of the target function. 

\section{A smooth variant using the penultimate right singular vector} 
\subsection{Motivating quotations}
From \cite{nakatsukasa2018aaa}, in 2018:
\begin{quote}
``If $|x|$ is approximated on a set $Z \subseteq [-1,2]$, for example, then all support points will of course be real and no complex numbers will appear. Nevertheless, there is the problem that the AAA algorithm increases the rational type one at a time, so that odd steps of the iteration, at least, will necessarily have at least one real pole, which will often appear near the singular point $x = 0$. We do not have a solution to offer to this problem.'' 
\end{quote}
From \cite{huybrechs2023aaa}, in 2023:
\begin{quote}
``It does seem that further research about avoiding unwanted poles in AAA approximation is called for, but fortunately, in a practical setting, such poles are immediately detectable and thus pose no risk of inaccuracy without warning to the user.'' 
\end{quote}
From \cite{driscoll2024aaa}, in 2024:
\begin{quote}
``Finally there is the perennial question in AAA approximation — indeed, in rational approximation generally — of what to do when there are poles in a region where $r$ should be analytic... Our experiments show that in difficult cases, bad poles often arise in approximately alternate steps...'' 
\end{quote}

\subsection{The AAAsmooth variant}
In this variant, the inputs from the user are the same as in AAA. The initialization and greedy addition of support points are also unchanged. 
To compute the weights at each step, construct the Loewner matrix as above, using data at $N$ support points $z_j$ and $M$ sample points $Z_i$. Compute its reduced SVD $A = U \Sigma V^H$, so that $\Sigma$ and $V$ are square $N\times N$ matrices; we do not need the $M\times N$ factor $U$. Let $\sigma_{N-1}$ and $\sigma_N$ be the two smallest singular values, and let $V_{N-1}$ and $V_N$ be the corresponding right singular vectors, i.e. the last two columns of $V$. Now instead of taking $w=V_N$, as in standard AAA, we use the following linear combination of $V_{N-1}$ and $V_N$:  
\begin{equation}
w = V_N + \left(\sigma_N/\sigma_{N-1}\right)^{\kappa} \sqrt{-1}\,  V_{N-1},
\label{eq:wi}
\end{equation}
where $\kappa=3/2$. 
\subsection{Discussion} 
The AAA algorithm described in the Introduction carries a special hazard for real-valued applications: it can produce poles on the real line, possibly within the domain of approximation. And `can' becomes `must' when there are an odd number of poles since they appear in complex-conjugate pairs for real problems. Every second step of the AAA process is especially vulnerable to the appearance of spurious poles, which is a problem not just when $f$ is smooth but also when $f$ has a finite number of poles (due to parity mismatch with the pole count of $r$). One strategy that ameliorates this issue is to make the real problem complex from the outset, say by approximating $f(x) + i e^x$ instead of $f(x)$, and then using only the real part of the approximation. 

The AAAsmooth variant described here originated as another strategy for circumventing this hazard in real cases. It also detours into the complex plane, although we will be able to return a real rational function for real problems if desired. A pleasant surprise is that this variant also improves on the AAA process for fully complex problems, so we continue to use complex notation, e.g. writing $f(z)$ instead of $f(x)$. 

The equation \eqref{eq:wi} deserves some discussion.  The recommendation $\kappa=3/2$ is a heuristic; the alternatives $\kappa=2$,  $\kappa=1$, and $\kappa = 1/2$ lead to fairly similar behavior and even $\kappa=0$ is only slightly worse (e.g. for the tests in Figure \ref{fig:twelve} below). The motivation for the expression $(\sigma_N/\sigma_{N-1})^\kappa$ is that the usual AAA choice $w = V_N$ is expected to be no better than $w = V_{n-1}$ when $\sigma_{N-1}\approx\sigma_N$, in which case we combine both versions with roughly equal weight; on the other hand, $w=V_N$ is superior if $\sigma_N \ll \sigma_{N-1}$. Thus \eqref{eq:wi} returns an answer very similar to that of standard AAA if $\sigma_{N}$ and $\sigma_{N-1}$ are well separated. In fact the original idea that led to this variant was to take $w = \theta V_N + (1-\theta)V_{N-1}$ and optimize the choice of $\theta\in[0,1]$. However, this often led to adjacent nodes having weights of the same sign, which guarantees that a real pole will appear. To circumvent this we then tried interpolation between $V_N$ and $\sqrt{-1}V_{N-1}$ (since scaling of barycentric weights by a constant has no effect, $w=V_{N-1}$ and $w=\sqrt{-1}V_{N-1}$ give the same rational approximant), which eventually led to \eqref{eq:wi}. 

For real cases, this trick complexifies the problem so that poles need not appear in complex conjugate pairs. In particular, there can be an odd number of poles without any requirement that one of them appear on the real line. On the other hand, the complexification is confined to the weights; the function values $f_j$ and nodes $z_j$ are still real. Therefore, if we wish to evaluate $r(z)$, defined by \eqref{eq:r(z)}, for $z=x\in \mathbb{R}$, we can express its real part as a rational function of $x$: 
\begin{equation}
\Re[r(x)] 
= 
\frac{
\left(\displaystyle\sum_j\frac{f_j(V_N)_j }{x-z_j}\right)
\left(\displaystyle\sum_j\frac{(V_N)_j }{x-z_j}\right)
+
\left(\displaystyle\frac{\sigma_N}{\sigma_{N-1}}\right)^{2\kappa}\left(\displaystyle\sum_j\frac{f_j (V_{N-1})_j}{x-z_j}\right)
\left(\displaystyle\sum_j\frac{  (V_{N-1})_j }{x-z_j}\right)
}{\left(\displaystyle\sum_j\frac{(V_N)_j }{x-z_j}\right)^2
+\left(\displaystyle\frac{\sigma_N}{\sigma_{N-1}}\right)^{2\kappa}\left(\displaystyle\sum_j\frac{  (V_{N-1})_j }{x-z_j}\right)^2}
\label{eq:realpart} 
\end{equation}
This can be returned as the purely real answer to the original, real rational approximation problem; note, however, that it is of higher degree than $r$ and agrees with the real part of $r$ only on the real line. For applications where complex arithmetic is acceptable it is probably better to use $r(z)$ or even $\Re[r(z)]$ as the approximant; in the numerical results presented here we always report errors with imaginary part included, that is, we report $|r-f|$ rather than $|\Re(r)-f|$. In the formula \eqref{eq:realpart}, it is also clear that the denominator can vanish only if both $\sum_j\frac{(V_N)_j }{x-z_j}$ and $\sum_j\frac{(V_{N-1})_j }{x-z_j}$ vanish. In contrast, the usual AAA approximant has the denominator $\sum_j\frac{(V_N)_j }{x-z_j}$. Thus it is more difficult, although not strictly impossible, for the AAAsmooth process to produce a real pole. This did not occur at all in our testing (although poles of very small imaginary part do appear when $f$ has a real pole).  

In problems that are complex-valued from the outset, the behavior of this variant is broadly similar to the AAA process but it is often much less volatile. By this we mean that the error curves for this variant, formed by plotting the maximum of $|f-r|$ on the sample grid at each step, are monotone decreasing or nearly so, in contrast to those of the standard AAA which can exhibit noisy behavior. For this reason we call the variant `{\tt AAAsmooth}'; the name also evokes the original motivation of avoiding poles in real problems with smooth target functions. The explanation for this behavior is unknown; one happy result of the less noisy error curve is that the AAAsmooth process often reaches the error tolerance one or several steps earlier than standard AAA, which may reduce the appearance of Froissart doublets. 

It is reasonable to ask if the AAAsmooth process is suitable for situations where $f$ is known to have poles on the real line, since it is much more difficult for this variant to produce real poles. The answer is affirmative if we accept approximations whose poles have small imaginary parts. As an example, we reexamine the very first application of the AAA process in \cite{nakatsukasa2018aaa}, the approximation of $\Gamma(z)$ from 100 evenly spaced samples in $[-1.5,1.5]$ with a tolerance of $10^{-13}$. The poles of the approximants produced by both AAA and AAAsmooth at steps 3, 4, 5, and 7 are listed in Table \ref{tbl:gammapoles}. The AAA process produces spurious poles at steps 3 and 5, but not thereafter (both versions terminate with 11 poles, or $N=12$). In the AAAsmooth process, the ratio $\sigma_{N-1}/\sigma_N$ takes the successive values $1.01$,
$2260$,
$6.22$,
$148$,
$4.42$,
$268$,
$150$,
$42.4$,
$213$, and 
$175$; when that ratio is high, AAAsmooth and AAA would select nearly the same barycentric weights. Thus there are only a few steps where AAAsmooth deviates significantly from AAA. 
The AAAsmooth process can yield poles that are nearly on the real line, e.g. $-2.0000333+0.0000007i$ at step 7. The imaginary part is usually smaller than real part of the error in these pole locations. Overall, there is little difference between the quality of the approximants produced by the two methods at the termination step, although at some intermediate steps the AAAsmooth approximant is better. 

\begin{figure}
\[\includegraphics{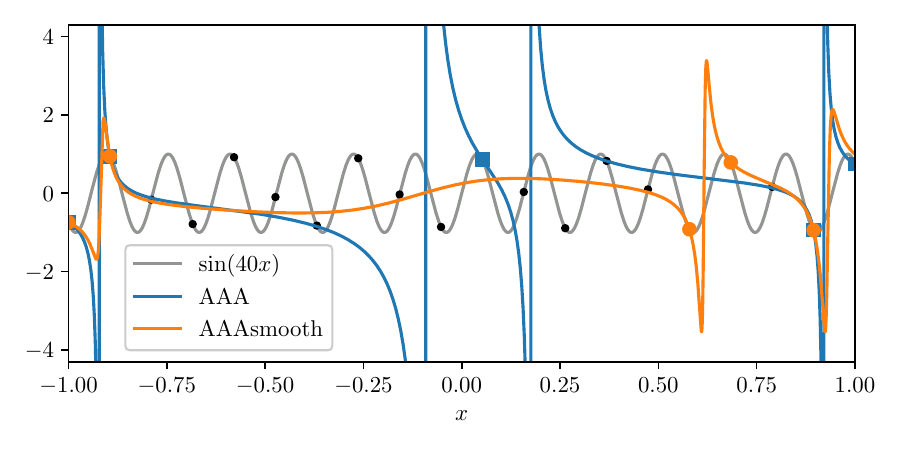}\]
\caption{Two barycentric rational approximations of $f(x)=\sin(40 x)$, each with five support points. The data set consists of 20 equispaced points on $[-1,1]$; the support points chosen by AAA are marked with blue squares while those chosen by AAAsmooth are marked with orange circles. Black dots indicate the remaining sample points. The AAA version has four poles within the interval while the AAAsmooth version has four complex-valued poles which lie within a distance of $0.01$ from the interval. }
\label{fig:stepfive}
\end{figure}

It is possible for AAAsmooth to produce poles that are close to the real line even when this is undesirable. Figure \ref{fig:stepfive} shows the AAA and AAAsmooth approximants that are produced for $f(x)=\sin(40x)$ using twenty evenly-spaced samples on $[-1,1]$, at the stage when five support points have been chosen. We only plot the real part of the AAAsmooth approximant. The AAA approximant has four real poles on $[-1,1]$, the worst possible case. As the figure suggests, the poles of the AAAsmooth variant are close to $[-1,1]$, specifically near $-0.917-0.009i$,  $-0.105-0.320i$,  $0.930+0.009i$, and 
  $0.616-0.006i$. While the AAAsmooth curve remains finite and is preferable to the AAA curve, this figure indicates that the AAAsmooth variant does not offer a completely satisfying answer to the problem of real-valued poles in AAA approximation. 
  
\begin{table}
\begin{tabular}{m{2.5cm}m{3.0cm}m{3.5cm}m{3.5cm}m{3.5cm}}
\toprule
Variant & Step 3 & Step 4 & Step 5 & Step 7\\
\midrule
AAA & 
$-0.999978$\par \vspace{0.12cm}
$\textcolor{blue}{-0.817512}$\par \vspace{0.12cm}
$-0.000080$\par \vspace{0.12cm}
&
$-1.9901991984$\par \vspace{0.12cm}
$-0.9999999969$\par \vspace{0.12cm}
$-0.0000001102$\par \vspace{0.12cm}
$3.2934195661$\par \vspace{0.12cm}
&
$-1.9885705533$\par \vspace{0.12cm}
$-0.9999999989$\par \vspace{0.12cm}
$0.0000001785$\par \vspace{0.12cm}
$\textcolor{blue}{0.6647336258}$\par \vspace{0.12cm}
$3.3442023766$\par \vspace{0.12cm}
&
$-2.94476117441128$\par \vspace{0.12cm}
$-2.00005523139197$\par \vspace{0.12cm}
$-0.99999999999987$\par \vspace{0.12cm}
$0.00000000000018$\par \vspace{0.12cm}
$3.98275290176057$\par$\pm2.08720414072944i$\par\vspace{0.12cm}
$4.32084281204053$\par \vspace{0.12cm}
\\
\midrule
AAAsmooth
&
$-1.76130361$\par$+0.09487590i$\par\vspace{0.12cm}
$-0.99998851$\par$-0.00000300i$\par\vspace{0.12cm}
$-0.00000072$\par$-0.00000017i$\par\vspace{0.12cm}
&
$-1.99373314820$\par$-0.00045363470i$\par\vspace{0.12cm}
$-1.00000020399$\par$-0.00000001522i$\par\vspace{0.12cm}
$0.00000000008$\par$+0.00000000086i$\par\vspace{0.12cm}
$3.32574427836$\par$+0.00527764307i$\par\vspace{0.12cm}
&
$-1.9874544755$\par$+0.0005713988i$\par\vspace{0.12cm}
$-0.9999999973$\par$+0.0000000014i$\par\vspace{0.12cm}
$-0.0000000063$\par$+0.0000000017i$\par\vspace{0.12cm}
$0.9459901704$\par$+0.5165947798i$\par\vspace{0.12cm}
$3.2708990563$\par$-0.0290823504i$\par\vspace{0.12cm}
&
$-2.95075213620777$\par$-0.00031789042078i$\par\vspace{0.12cm}
$-2.00003330798366$\par$+0.00000070511062i$\par\vspace{0.12cm}
$-1.00000000000043$\par$+0.00000000000001i$\par\vspace{0.12cm}
$0.00000000000043$\par$-0.00000000000002i$\par\vspace{0.12cm}
$3.90953039063703$\par$-2.05993708767955i$\par\vspace{0.12cm}
$3.91767605794173$\par$+2.03807520476246i$\par\vspace{0.12cm}
$4.24630330769600$\par$-0.01202478749298i$\par\vspace{0.12cm}
\\
\bottomrule
\end{tabular}
\caption{The poles of the successive rational approximations generated by the AAA and AAAsmooth algorithms while approximating $\Gamma(z)$ using data at 100 equally spaced values on $[-1.5,1.5]$. Both versions identify the poles (at 0 and negative integer values) to increasing accuracy as the algorithm progresses; by step 7 both have found the poles at $-2$ and $-3$ even though these lie outside the data interval. 
The AAA approximation produces spurious poles within $[-1.5,1.5]$ at steps 3 and 5 (blue text) whereas the AAAsmooth process does not. However, there are no spurious poles after step 5; in particular, at step 7 the two versions give very similar results. Note that in Step 7, the AAAsmooth poles nearly have conjugate symmetry, since at that stage we have $\sigma_{N-1}/\sigma_N = 150$ which leads to nearly real barycentric weights.  }
\label{tbl:gammapoles}
\end{table}

We present more numerical examples later, after introducing a second AAA variant. 

\section{A cheap variant using first derivative information}
\subsection{Motivating quotation}
The following quotation from \cite{nakatsukasa2018aaa}, published in 2018, in turn refers to an older work \cite{antoulas1986scalar}: 
\begin{quote}
``\emph{Confluent sample points} is the phrase we use to refer to a situation in which one
wishes to match certain derivative data as well as function values; one could also
speak of ``Hermite-AAA” as opposed to ``Lagrange-AAA” approximation. Antoulas
and Anderson consider such problems... [their setting] is interpolatory, so weighting of rows is not an issue except perhaps for numerical stability, but such a modification for us would require a decision about how to weight derivative conditions relative to others.''
\end{quote}
\subsection{The AAAbudget variant} 
In this variant the user must supply more information than in standard AAA: the values and first derivatives of $f$ are given on a finite sample set. The conversion of sample points to support points is the same as in AAA, both for initialization and for the greedy selection of new support points; those decisions do not use the derivative information. 
The barycentric weights are determined using the given information at the support points, ignoring the data at sample points. Specifically, given support points $z_1$, $\cdots$, $z_N$, we form the $N\times N$ matrix $B$ whose off-diagonal entries are the difference quotients $(f(z_i)-f(z_j)) / (z_i-z_j)$, and whose diagonal contains the derivatives $f'(z_i)$. We then compute the SVD of $B$ and take the entries of the last right singular vector as barycentric weights. 

\subsection{Discussion}
The formula for $r'(z_i)$ given above in \eqref{eq:r'(z)} can be rearranged to give $ 0=w_i r'(z_i) + \sum_{j\neq i} \frac{f_i-f_j}{z_i-z_j}w_j$. It is striking that most of the coefficients are first-order finite difference approximations of $f'$, like the Loewner matrix entries. 
When we enforce $r'(z_i)=f'(z_i)$ for $i=1,\cdots, N$ we obtain the equations 
\begin{equation}
0=w_i f'(z_i) + \sum_{j\neq i} \frac{f_i-f_j}{z_i-z_j}w_j
\label{eq:r'cond}
\end{equation}
which simply state that $w$ lies in the nullspace of the $N\times N$ matrix $B$ defined above. Thus, we are seeking weights that will give us desired first derivative information at support points, rather than desired function values at sample points. In this admittedly narrow setting we find that we can sidestep the reweighting question mentioned in the quotation above. That is, for the very specific case where we use first-order derivative information at the support points, the entries of $B$ all have the same dimensions (which are the same dimensions as the usual Loewner matrix entries) and no row reweighting is needed. 

We now have a natural way to include first-order derivative information in a AAA-style algorithm. But what do we gain by doing so? Initially we hoped that the use of this additional information would allow us to improve on the AAA process. We tried several variants in which we combined the information in the $N\times N$ matrix $B$, which minimizes $|r'-f'|$ at support points, with the $M\times N$ Loewner matrix $A$, which minimizes $|r-f|$ at sample points, e.g. stacking these matrices before taking the SVD. None of these versions performed better than the standard AAA process, except in cases where relatively few sample points were given; in this case the convergence of AAA was limited by the available data, whereas a competing variant with access to twice as much information (values of $f'$ as well as $f$ at all samples) could do better. This is not a fair comparison; we address this issue more precisely below with Figure \ref{fig:even}. We eventually settled on the version above because it has the advantage of being much faster than AAA when there are many sample points, while giving similar results. The name `AAAbudget' highlights the decreased cost of this variant. Indeed most of the computational expense of the algorithm comes from the SVD, and most of the SVD time is incurred in the last few steps of the algorithm when the matrices are large. If the cost of an SVD is $\mathcal{O}(MN^2)$ for AAA versus $\mathcal{O}(N^3)$ for AAAbudget, and both versions conclude with $N$ support points, and $M\gg N\gg 1$, we expect AAAbudget to be faster by a factor of $\mathcal{O}(M/N)$. Of course, the AAA procedure is usually very fast already, frequently terminating in under a second; moreover, the continuous version of AAA described in \cite{driscoll2024aaa}, with adaptive sampling, is faster than standard AAA in many situations. Thus, the new AAAbudget variant adds the most value, relative to existing methods, in the very special case where where the approximation must hold on a region rather than a curve in $\mathbb{C}$, AAA must be called repeatedly, and first derivative information is available. 

Before we present numerical examples we comment further on equation \eqref{eq:r'cond}. 
We might hope, despite the nonlinear dependence of $r'$ on $w$, that we can always choose the $N$ barycentric weights to satisfy the $N$ first-order conditions $r'(z_j)=f'(z_j)$. In general this hope is unfounded. Indeed we have already observed that if the data values are all identical, then $r(z)$ is a constant function and the choice of the barycentric weights has no effect; thus there is no way to satisfy any derivative conditions other than $r'(z_j)\equiv0$. As a second example, consider the problem of approximating $f(z) = |z| = z^2/\sqrt{z^2}$ by a barycentric rational function with the four support points $\{-3,-1,1,3\}$. Then the matrix $B$ has the simple form 
\begin{equation}
B = \begin{pmatrix}
-1&-1&-\frac12&0\\
-1&-1&0&\frac12\\
-\frac12&0&1&1\\
0&\frac12&1&1\\
\end{pmatrix}.
\end{equation}
Here $B$ is nonsingular and there is no way to satisfy all four first derivative conditions. Our procedure imitates AAA by adopting the entries of the last right singular vector as barycentric weights; an alternative that we do not examine here would be to find the smallest possible perturbation of the diagonal of $B$ such that it becomes singular, and then use the nullspace vector of the perturbed matrix as the barycentric weights. A third, contrasting example is the problem of approximating $f(z) = \sqrt{1.21-z^2}$ at the five points $\{0,\pm0.5,\pm1\}$. In this case the matrix $B$ has a one-dimensional nullspace, leading to an approximant $r(z)$ whose first derivatives perfectly match those of $f$ at those five points. The error curve $r-f$ therefore vanishes and also has a vanishing first derivative at the support points; see Fig. \ref{fig:surprise1}. This results in a curious behavior which is in some sense the opposite of the `equioscillation' displayed by rational best approximants: $r(z)$ is an \emph{upper bound} for $f$, that is, $r(z)\ge f(z)$ for all $z\in[-1,1]$, with equality only at the support points. Thus the matrix $B$ can be nonsingular or singular depending on which derivative conditions we are seeking to satisfy. From a linear algebra viewpoint this makes sense: any square matrix, including our $B$, may be made singular or nonsingular if we are free to change its diagonal entries. The corresponding statement about $r$ is that it may or may not be possible to choose the weights to satisfy $r'(z_j)=f'(z_j)$, depending on the values $f'(z_j)$ that are prescribed. If it is impossible, we can still hope to get nearly the correct derivative values by using the entries of the last right singular vector for the barycentric weights.\footnote{One can also use a complex linear combination of the last two right singular vectors of $B$, thereby combining the AAAsmooth and AAAbudget ideas. However, the resulting algorithm does not meaningfully improve on AAAbudget.} 
\begin{figure}
\[\includegraphics{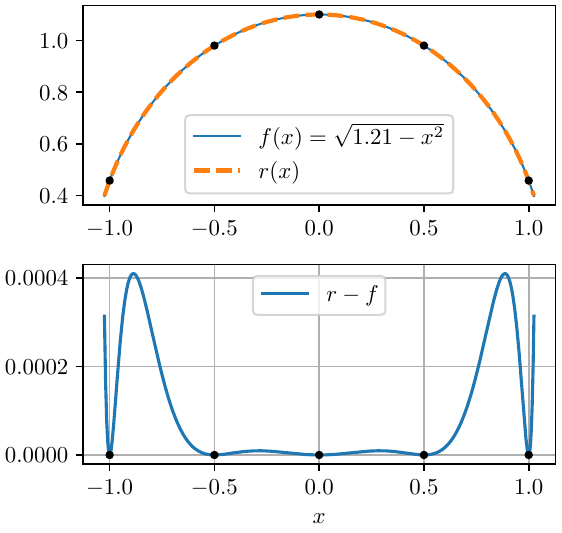}\]
\caption{This approximation $r(x)$ 
of $f(x)=\sqrt{1.21-x^2}$ is also an upper bound for $f$. The values and first derivatives of $r$ and $f$ agree at the five support points $\{-1,-0.5,0,0.5,1\}$. The top panel shows $r$ and $f$, while the lower panel shows the difference $r-f$ (which is nonnegative). }
\label{fig:surprise1}
\end{figure}

We close this section by noting that a zero barycentric weight is less potentially catastrophic for the AAAbudget variant than for the AAA or AAAsmooth variants. Indeed, suppose that \eqref{eq:r'cond} holds exactly and it happens that $w_i=0$. The resulting statement, 
\begin{equation}
0=\sum_{j\neq i} \frac{f_i-f_j}{z_i-z_j}w_j
\end{equation}
is precisely the equation for ensuring that $n(z_i)=f_i d(z_i)$, where $n$ and $d$ are defined without the term involving $z_i$. In other words, if $w_i=0$ then we lose control of $r'(z_i)$ but we retain the zeroth-order condition $r(z_i)=f_i$. Even if \eqref{eq:r'cond} only holds approximately, we can expect $|r(z_i)-f_i|$ to be no larger than the errors in $r'$ at other points.  

\section{Numerical Tests}
We now compare the performance of the standard AAA approximation with the AAAsmooth and AAAbudget variants, focusing on some examples that have already been considered in other recent works. The computations are performed in Python; our implementation of AAA can reproduce fine details of some figures created by other authors with the Chebfun version (see Figure \ref{fig:even} below). For many problems, the three methods give similar results. We take a close look at three cases where they differ more than usual and then give a briefer treatment of twelve problems on which their performance is similar. 

\begin{figure}
\[\includegraphics{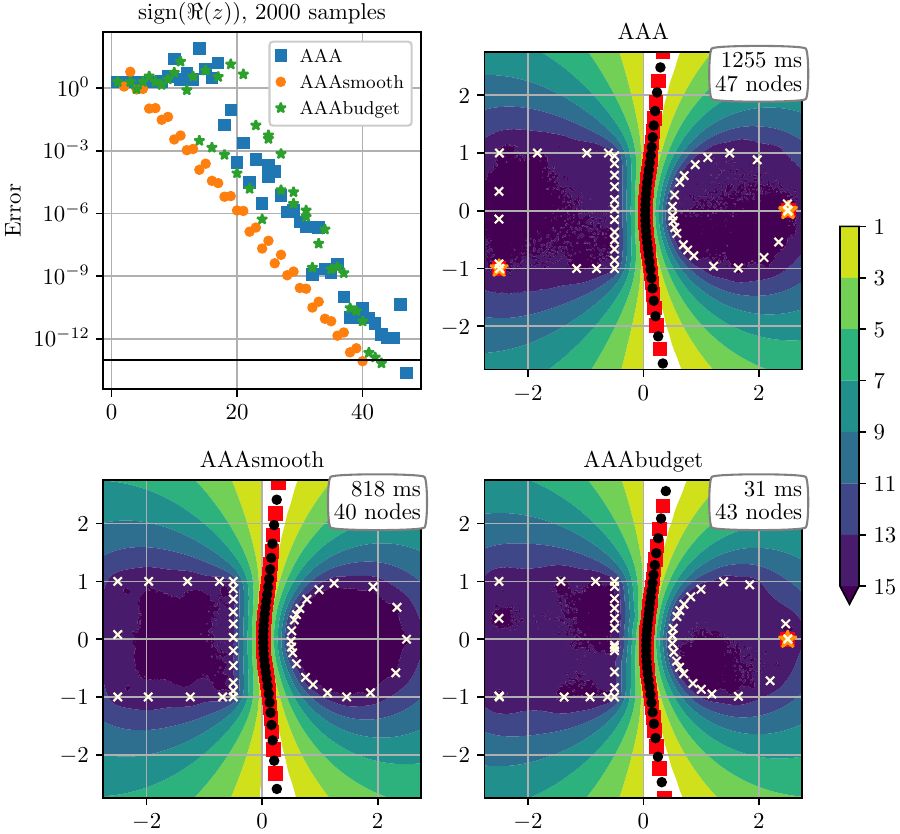}\]
\caption{A new look at Application 6.6 from \cite{nakatsukasa2018aaa}. We take $f=-1$ at 1000 points evenly spaced around a square on the left and $f=1$ at 1000 points evenly spaced around a circle on the right. The AAA algorithm constructs a rational approximation $r$ with 47 support points, marked with white `x' symbols in the upper right panel. The contour plots depict the number of correct decimal digits, $-\log_{10}|r-f|$; white zones indicate $|r-f|>0.1$. The zeros (red squares) and poles (black circles) of $r$ appear on a contour separating the square and the circle. The rational approximations produced by the AAAsmooth and AAAbudget variants are similar, but both have slightly fewer nodes than the AAA output. Orange stars mark the locations of numerical Froissart doublets (detected wherever a pole of $r$ of residue less than $10^{-10}$ times the largest residue appears inside a circle centered at a zero of $r$ with radius of $10^{-10}$ times the diameter of the support set). The AAAbudget process is 40 times faster than traditional AAA because of the large sample size (31 milliseconds versus 1.26 seconds). In the upper left panel the convergence curves for each variant indicate the largest value of $|r(z_k)-f_k|$ as a function of the number of support points. The AAAsmooth curve, while not quite monotone, is less noisy than the others and is the first to reach the tolerance of $10^{-13}$ on the set of sample points. }
\label{fig:square2circ}
\end{figure}

We begin with a piecewise constant function with a discontinuity along the imaginary axis: $f(z) = \sgn(\Re(z))$. The algorithm sees only that $f\equiv -1$ on a square and $f\equiv 1$ on a circle; the square has side length 2 and center $-1.5$ and the circle has radius $1$ and center $1.5$. Both the square and the circle are represented by 1000 evenly spaced points. This problem appeared as Application 6.6 in \cite{nakatsukasa2018aaa}. The error curves (plots of $\max_i|f(Z_i)-r(Z_i)|$ as a function of the number of support points) are given in the upper right panel of Fig. \ref{fig:square2circ}; the other panels illustrate the final version of $r$ returned at convergence by the three algorithms. In all three cases the approximant $r$ is nearly constant on the interiors of the square and circle, and the line of discontinuity is approximated by an alternating sequence of zeros and poles of $r$. The contours indicate agreement of $r$ with $f$, with each contour line representing two correct decimal digits; the meaning of the colors is consistent across all of the contour plots in Figures \ref{fig:square2circ}, \ref{fig:expessential}, and \ref{fig:tan256}. The error curves have similar slopes but the AAAsmooth curve is much less noisy than the others, and it is the first to reach the tolerance of $10^{-13}$. Thus the AAAsmooth variant gives the most compact barycentric representation among the three and it is the only one where Froissart doublets do not appear. Code timings are also reported for the three versions; these times are the averages of 40 runs as measured with Python's {\tt time.process\_time} routine. These timings do not include the computation of poles and zeros and residues, although that computation is not a large component of the overall cost. Because the number of sample points is large, the AAAbudget version has a significant speed advantage, requiring only about 1/40th of the time needed for standard AAA. The AAAsmooth process is also somewhat faster than AAA because it converges in fewer iterations, thereby avoiding the most expensive SVD calls. 
This example stands out from those in the later Fig. \ref{fig:twelve} in two ways. It is not surprising that AAAsmooth terminates first, but often the difference between AAA and AAAsmooth is 0-2 iterations rather than 7. Also, it is unusual that the AAAbudget procedure takes fewer steps than AAA to converge; usually AAAbudget requires a few more iterations to converge. 

\begin{figure}
\[\includegraphics{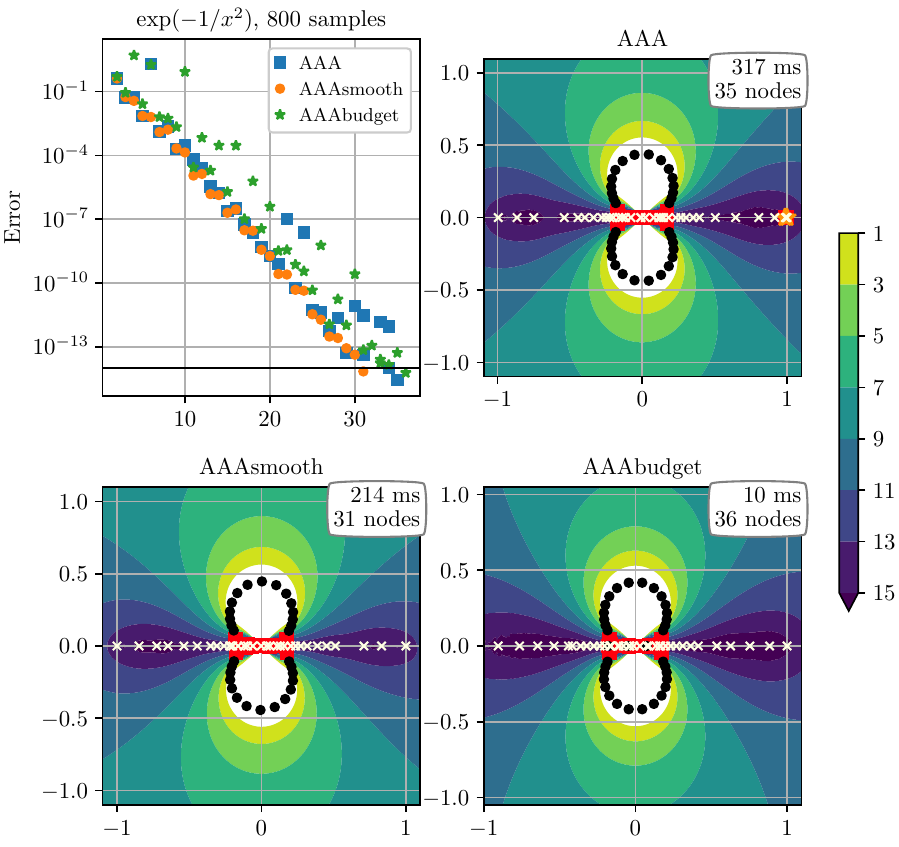}\]
\caption{A comparison of AAA, AAAsmooth, and AAAbudget for approximating $f(z)=\exp(-1/x^2)$ in 800 evenly spaced points on $[-1,1]$ with tolerance $\epsilon=10^{-14}$. The white regions in the contour plots indicate regions where $|r-f|>1$. Note that AAAsmooth has a less noisy error curve than the other methods. }
\label{fig:expessential}
\end{figure}

Our second example is a real-valued problem, the approximation of the function 
$f(z) = \exp(-1/z^2)$ 
on $[-1,1]$ in 800 evenly spaced points with tolerance $\epsilon=10^{-14}$; see Fig. \ref{fig:expessential}. If we restrict the domain of definition to $\mathbb{R}$, this function is infinitely differentiable but not analytic at the origin; if the domain is $\mathbb{C}$ we would instead say that it has an essential singularity at $z=0$. It appears as an example in \cite{driscoll2024aaa} and \cite{huybrechs2023aaa}. Again, the error curves for AAA and AAAbudget are noisy compared to the curve for AAAsmooth, which converges with 31 support points instead of 35 (AAA) or 36 (AAAbudget). With 800 samples, AAAsmooth has a speed advantage, requiring only 10 ms compared to 317 ms for AAA and 214 ms for AAAsmooth. The white regions in the contour plots indicate regions where $|r-f|>10^{-1}$. From the error curves we can see that the choice $\epsilon=10^{-13}$ would have caused both AAA and AAAsmooth to terminate with 29 nodes; despite the noise in the AAA convergence curve, many of its iterations have errors as small as those of AAAsmooth.  

\begin{figure}
\[\includegraphics{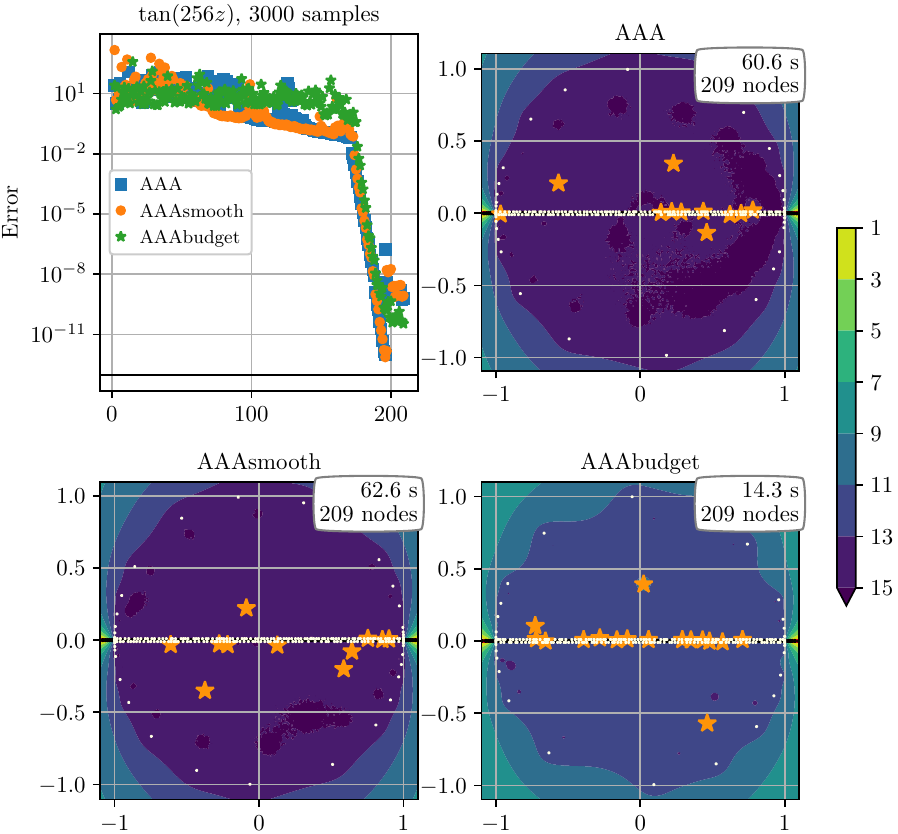}\]
\caption{Three approximations of $f(z)=\tan(256z)$ in 3000 points, of which 1000 are evenly spaced on the unit circle and 2000 are evenly spaced on the two horizontal segments with endpoints $\pm1 \pm0.01i$. None of these methods achieved the convergence criteria and the iterations all continue until the imposed maximum of 210 support points (white dots). This stagnation is associated with the appearance of Froissart doublets (stars). While the AAA and AAAsmooth results have better accuracy than AAAsmooth on most of the unit disk, AAAsmooth actually has a smaller maximum error overall, as the convergence curve shows. }
\label{fig:tan256}
\end{figure}

In the third example we consider a problem where AAA and the two variants all fail to reach the tolerance of $10^{-13}$. The task is to approximate $f(z)=\tan(256z)$ on the unit disk, based on data at 1000 equally spaced points on the boundary of the disk as well as 2000 more points which are equally spaced on two horizontal line segments joining $\pm 1\pm 0.01i$. This is similar to Application 6.4 in \cite{nakatsukasa2018aaa}, but with the two segments replacing randomly distributed points in the disk. The (absolute) error tolerance $\epsilon=10^{-13}$ is somewhat too tight because $|f|$ has a maximum value of about 25 on this 3000-point grid: we are asking for more than 14 correct digits, which is at the outer edge of what the AAA process can achieve in double precision arithmetic. As the top left panel of Fig. \ref{fig:tan256} shows, all three versions show good convergence up to almost $N=200$ but then the error curves jump up by several orders of magnitude; we artificially halt the iteration at $N=210$. The contour plots state that $N=209$, not $N=210$, because one node had a vanishing weight for each variant. This is a rare example where the vanishing weight is dangerous. For the AAA process the vanishing weight at the last step corresponds to the first node, at $z=1$, which is very close to a pole at $z=163\pi/512=1.000155$. Without the interpolative property to guarantee $r(1)=f(1)$, we find that $|r(1)-f(1)|=6.1\cdot10^{-10}$ whereas the maximum of $|r-f|$ at the other 2999 points is $1.3\cdot10^{-13}$. That is, the convergence criterion was nearly satisfied except at the one support point whose weight vanished. The picture is very similar for the AAAsmooth variant, but in line with our comments above the AAAbudget process is less vulnerable to the hazard of vanishing weights. With the AAAbudget approximation, we find an error of $6.5\cdot10^{-12}$ at the node whose weight vanishes and errors up to $3.5\cdot10^{-11}$ on the rest of the grid: this is a better outcome in the $\infty$-norm, although the contour plots indicate that AAAbudget gives a worse result on most of the disk.

\begin{figure}
\[\includegraphics{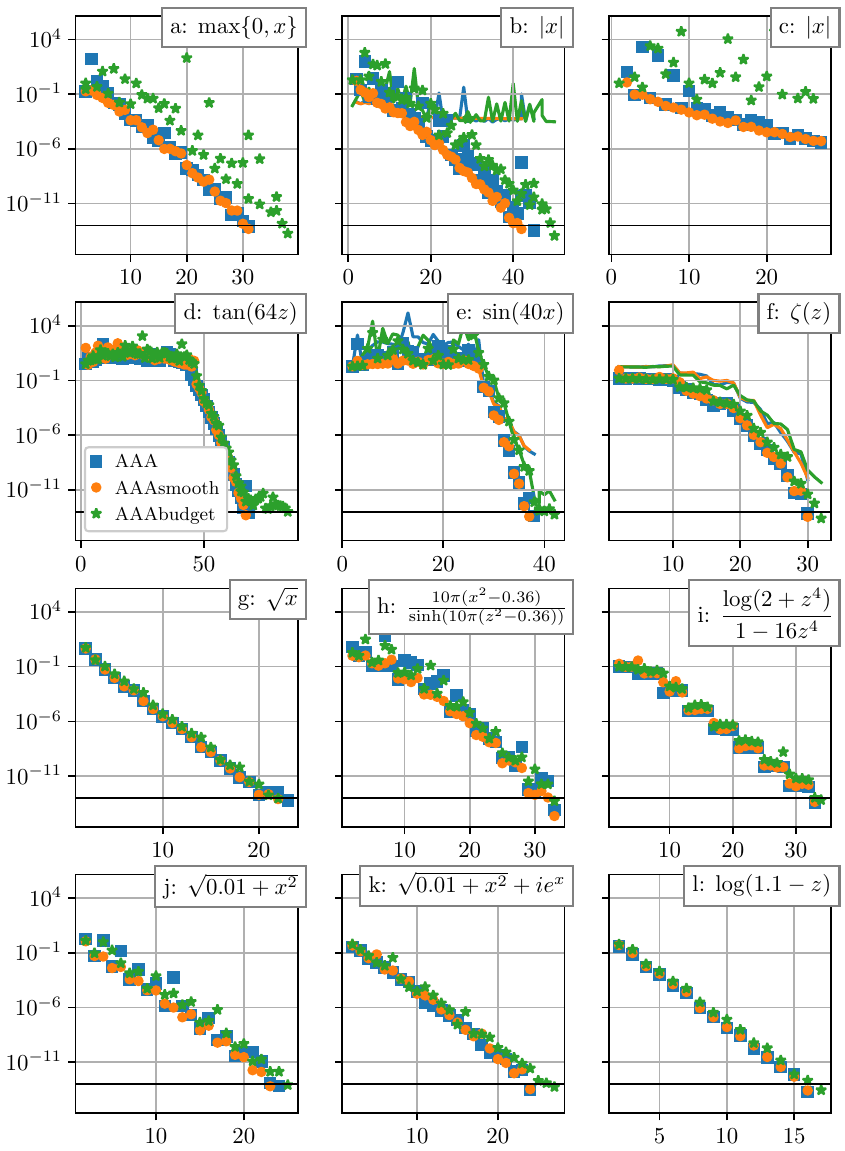}\]
\caption{The AAA, AAAsmooth and AAAbudget versions often have similar behavior; see text for details. }
\label{fig:twelve}
\end{figure}

Finally in Fig. \ref{fig:twelve} we give convergence curves for twelve other problems which also appear in the recent literature on the AAA algorithm. In all cases the vertical axis shows $|r-f|$ and the horizontal axis gives the number of support points, which is typically two more than the iteration count since the algorithm begins with two support points. We write $f(x)$ instead of $f(z)$ if all of the data samples are taken on the real line. The tolerance is $\epsilon=10^{-13}$ in all cases. More information on each example follows: 
\begin{itemize}
\item[(a)] Approximation of $f(x) = \max\{0,x\}$ in 200 evenly spaced points on $[-1,1]$, following \cite{driscoll2024aaa}. 
\item[(b)] Approximation of $f(x) = |x|$ in 1001 evenly spaced points on $[-1,2]$, following a comment in Section 8 of \cite{nakatsukasa2018aaa}. This grid avoids 0 so that we do not need to give a value for $f'(0)$. In addition to the usual convergence curves given by discrete markers, we plot solid lines showing the maximum errors on a fine evaluation grid of 1000 points evenly spaced on $[-0.01,0.01]$; this shows that all three methods still make large errors near the origin (in between the given samples) when convergence is reached. 
\item[(c)] Approximation of $f(x) = |x|$ in 200000 evenly spaced points on $[-1,1]$, following Application 6.7 in \cite{nakatsukasa2018aaa}. Here we plot only the initial iterations, up to 27 support points, to imitate Figure 6.10 in \cite{nakatsukasa2018aaa}. While not shown, the AAA and AAAsmooth iterations later reach the tolerance of $10^{-13}$ after about 70 steps (the AAA curve is smooth from iteration 15-55 and then becomes noisy again for $55<N<70$) and the AAAbudget curve stagnates, with errors between $10^{-8}$ and $10^{-12}$ after iteration 70.  
\item[(d)] Approximation of $f(z) = \tan(64z)$ in 3000 points, of which 1000 are evenly spaced around the unit circle and the other 2000 are evenly spaced on the two horizontal segments joining $\pm1\pm0.01i$. This imitates Application 6.4 in \cite{nakatsukasa2018aaa}. The only change is that our two horizontal segments have replaced 2000 uniformly randomly distributed points within the disk, since we found that different trials with different random placements can lead to meaningfully different outcomes. The maximum value of $|f|$ is about 25, so convergence to the (absolute) tolerance $\epsilon=10^{-13}$ really means better than 14 digit accuracy. This is at the edge of what can be achieved with AAA. With our fixed grid, the AAA and AAAsmooth variants converge nicely while the AAAsmooth variant exhibits some stagnation near the end. 
\item[(e)] Approximation of $f(x) = \sin(40x)$ in 90 evenly spaced points on $[-1,1]$, following \cite{huybrechs2023aaa}. In addition to the usual convergence curves given by discrete markers, we plot solid lines showing the maximum errors on a finer grid of 2000 equispaced points in $[-1,1]$. These solid curves show that AAA and AAAsmooth reach the tolerance criterion on the coarse grid while still making larger errors ($\mathcal{O}(10^{-8})$) on the finer grid. This is a case where the given information has nearly been exhausted by these two methods: with 39 support points, there are only 51 sample points remaining, so the shape of the Loewner matrix is starting to become uncomfortably close to a square. The AAAbudget process does better on the fine grid because it has access to twice as much information on the coarse grid (values of $f$ and $f'$). 
\item[(f)] Approximation of $f(z)=\zeta(z)$ on 100 points, evenly spaced between $4-40i$ and $4+40i$, following the Conclusion of \cite{nakatsukasa2018aaa}. On this sample set, $\zeta$ is not far from being a constant function: it satisfies $|1-\zeta|<0.082$. However, this information is enough for all three versions of AAA to locate the pole at $z=1$ and many nontrivial zeros with good accuracy (e.g. 10 digits). In addition to the usual convergence curves given by discrete markers, we plot solid lines showing the quantity $1/|r(1)| + |r(0.5 + 14.134725141734694i)|$ at each step.  
\item[(g)] Approximation of $f(x)=\sqrt{x}$ on 2000 points, geometrically spaced in $[0.01,100]$. This is Zolotarev's third problem; see \cite{filip2018rational,akhiezer1990elements}. 
\item[(h)] Approximation of $f(x) = 10\pi(x^2-0.36) / \sinh\left(10\pi(x^2-0.36)\right)$ in 160 evenly spaced points on $[-1,1]$. This function is a less difficult version of the $f_4$ used in \cite{filip2018rational}, which in turn cites \cite{van2010computing}; in those works the factors of $10$ are replaced with factors of $100$ and the focus is on finding the \emph{best} rational approximations of a given degree. 
\item[(i)] Approximation of $f(z)=\log(2+z^4) / (1-16 z^4)$ in 1000 points, equally spaced on the unit circle, following Section 5 of \cite{nakatsukasa2018aaa}, which in turn cites \cite{gonnet2011robust}. All three error curves improve meaningfully on every fourth step, since $f$ has four branch cuts, each approximated by alternating poles and zeros. 
\item[(j)] Approximation of $f(x) = \sqrt{0.01+x^2}$ in 80 evenly spaced points on $[-1,1]$. This example is from \cite{huybrechs2023aaa}. 
\item[(k)] Approximation of $f(x) = \sqrt{0.01+x^2}+\sqrt{-1}\exp(x)$ in 80 evenly spaced points on $[-1,1]$. This example is a complexified version of example (j); all reported errors are complex magnitudes, that is, we report $|r-f|$, not $|\Re(r-f)|$. The complexification leads to a somewhat smoother convergence curve for the AAA process. 
\item[(l)] Approximation of $f(z)=\log(1.1-z)$ in 256 points, evenly spaced on the unit circle, following Application 6.2 in \cite{nakatsukasa2018aaa}. 
\end{itemize}

One takeaway from Fig. \ref{fig:twelve} is that the performance of the three variants is very similar on many of the problems. Generally speaking, the AAAsmooth variant has less noisy error curves than the AAA variant, which in turn has less noisy error curves than the AAAbudget version. The examples in the first row make a clear suggestion that AAAbudget should not be used for piecewise-defined functions, at least when the point of discontinuity lies within the region where sample data is given. (The AAAsmooth variant does quite well in the example of Fig. \ref{fig:square2circ}, where the discontinuity does not approach the square and circle where the data are given). Comparing panels (j) and (k) we see that the complexification essentially removes the noise from the AAA error curve, making it resemble the AAAsmooth error curve. 

\section{Conclusion}
The main finding in this paper is that the AAAsmooth process gives identical or somewhat better performance than the standard AAA process, sometimes with a smoother decrease in errors at successive steps and with termination several steps earlier. The AAAsmooth process provides a way to avoid spurious poles in real-valued problems. The price is the introduction of complex arithmetic (or the use of a real rational function of much higher degree).  The AAAbudget process, which uses information on the first derivative of the target function, can be significantly faster than AAA, but in some problems the result is not as good. The AAAbudget variant can produce approximants which are also upper or lower bounds for the target function, which is unusual in rational approximation generally, and it is less vulnerable than AAA to the potential dangers that arise when a barycentric weight vanishes. 

Looking ahead, this paper invites further work to use more of the information provided by the singular value decomposition of the Loewner matrix. Is there a good use for the antepenultimate right singular vector, $V_{N-2}$? Can we use the left singular vectors $U_N$ and $U_{N-1}$, which contain information about the pointwise errors $d(z_j)f_j - n(z_j)$, to further optimize the choice of barycentric weights? Is there a good way to use first derivative information at sample points, that is, to enforce constraints on $r'$ even when $r$ may not be correct yet? 

We close with Figure \ref{fig:even}, which compares AAA with AAAsmooth and AAAbudget in the low-information context. This figure should be compared with Figure 2 in \cite{huybrechs2023aaa}, which shows the superior convergence of AAA compared to five other algorithms. In the first panel we take $f(x) = \sqrt{1.21-x^2}$, as in \cite{huybrechs2023aaa}, while in the second panel we use the sum of the other functions treated in Figure 2 of \cite{huybrechs2023aaa}, namely $f(x) = \sqrt{0.01+x^2}+\tanh(5x)+\sin(40x)+\exp(-1/x^2)$. 
In contrast with our earlier figures, the horizontal axis shows the total number of samples in some relatively coarse grids where data is given; there are $n$  evenly spaced points on $[-1,1]$, with $n=8,12,16,\cdots,200$. The vertical axis shows the error at the conclusions of many AAA-style processes (not the error at each iteration of a single call to AAA). While $|r-f|$ always falls below the tolerance of $10^{-13}$ on the coarse grids, the displayed errors are the maximum of $|r-f|$ on a finer evaluation grid of 1000 evenly spaced points in $[-1,1]$. Dots are printed whenever the approximant has a pole in $[-1,1]$. We do not use any mitigation procedure when poles appear on the real line, which is a difference from \cite{huybrechs2023aaa}. The AAAbudget curve, in green, shows much faster convergence as measured on the fine evaluation grid; this is because it has access to twice as much information as the AAA and AAAsmooth versions. The black dashed line is the same as the AAAbudget curve, except stretched by a factor of $2$ in the horizontal direction. This shows that the AAAbudget variant is similar to the others once the information advantage is accounted for. It is notable that the AAAsmooth process does not produce any real poles; moreover it only produces a few poles that are very close to the real line. Out of the 98 calls to AAAsmooth for this figure, only three yield a pole with real part in $[-1,1]$ and imaginary part within distance $10^{-5}$ of the real line, in the worst case as close as $9.8\times10^{-10}$; on an additional 10 occasions the imaginary part is within distance $10^{-3}$ of the real line. As an additional comment about our Python implementation of AAA, we remark that the small zigzag in the AAA curve in the left panel of Fig. \ref{fig:even} at $n=48$ and $n=52$ perfectly matches the same feature in Figure 2(A) of \cite{huybrechs2023aaa}, which was generated in MATLAB with the {\tt chebfun} implementation.

\begin{figure}
\[\includegraphics{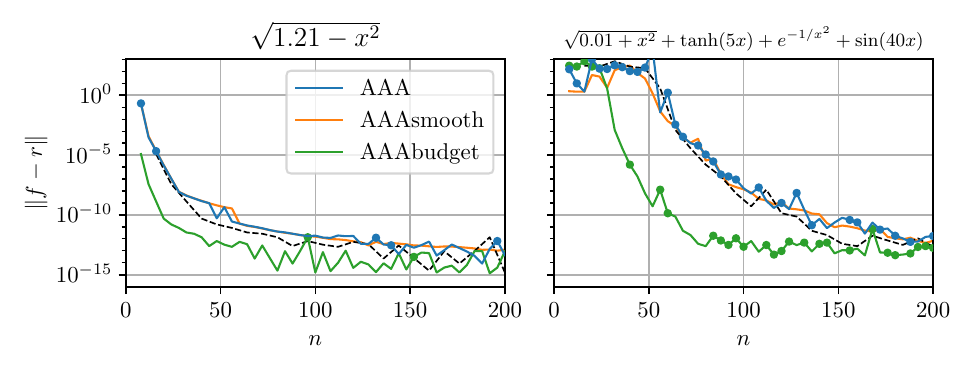}\]
\caption{Comparison of the three methods for approximation in evenly-spaced points, following Figure 2 in \cite{huybrechs2023aaa}. Here the horizontal axis indicates the number of points $n$ in an evenly-spaced grid where the data is given; the results are the outcomes of individual AAA calls, not progress towards convergence within a single call. Dots are printed wherever the approximant has a pole in $[-1,1]$. The black dashed line is the same as the green AAAbudget curve, except with $n$ doubled to make for a more fair comparison to the other methods since AAAbudget has access to $f'$ as well as $f$. 
The AAAsmooth variant never produces a real pole in $[-1,1]$. Also, the small zigzag in the AAA curve in the left panel near $n=50$, generated by our Python implementation, perfectly matches the same feature in \cite{huybrechs2023aaa}, generated in MATLAB with the {\tt chebfun} implementation. }
\label{fig:even}
\end{figure}

\bibliographystyle{plain}
\bibliography{mybib}

\end{document}